\documentclass[11pt,reqno]{amsart}
\usepackage{amsfonts}
\usepackage{amsmath,amsthm,amsfonts,amssymb}
\newtheorem{lemma}{Lemma}
\newtheorem{theorem}{Theorem}
\newtheorem{corollary}{Corollary}

\def\bl{\begin{lemma}}
\def\bt{\begin{theorem}}
\def\el{\end{lemma}}
\def\et{\end{theorem}}
\def\bp{\begin{proof}}
\def\ep{\end{proof}}
\def\bc{\begin{corollary}}
\def\ec{\end{corollary}}

\def\mb{\mathbb}

\def\O{\Omega}

\def\b{\beta}
\def\p{\partial}

\def\-{\setminus}

\def\vp{\varphi}
\def\ov{\overline}
\def\lt{\left}
\def\rt{\right}
\def\+{\bigcup}
\def\.{\bigcap}

\def\ll{\langle}
\def\rl{\rangle}

\title[A Schwarz-Pick lemma]{A Schwarz-Pick lemma for the modulus of holomorphic mappings between the
unit balls in complex spaces}
\thanks{* Corresponding author.\\
E-mail addresses: dymdsy@163.com (S. Dai), pan@ipfw.edu (Y. Pan)}

\begin{document}
\maketitle

\begin{center}
\bf Shaoyu Dai$^{a, *}$ , Yifei Pan$^{b, c}$

\bigbreak

\footnotesize\it  $^{a}$Department of Mathematics, Jinling Institute
of Technology, Nanjing 211169, China

\footnotesize\it  $^{b}$School of Mathematics and Informatics,
Jiangxi Normal University, Nanchang 330022, China

\footnotesize\it  $^{c}$Department of Mathematical Sciences, Indiana
University - Purdue University Fort Wayne, Fort Wayne, IN
46805-1499, USA

\end{center}

\numberwithin{equation}{section}


\hrulefill

\noindent {\bf Abstract:} In this paper we prove a
Schwarz-Pick lemma for the modulus of holomorphic mappings between
the unit balls in complex spaces. This extends the classical
Schwarz-Pick lemma and the related result proved by
Pavlovi$\acute{c}$.


\hrulefill

\section{Introduction}

Let $X$ be an open subset in the complex space $\mathbb{C}^n$ of
dimension $n$ and $Y$ be an open subset in the complex space
$\mathbb{C}^m$ of dimension $m$. Let $\mathbb{B}_n$ be the unit ball
in $\mathbb{C}^n$. The unit disk in the complex plane is denoted by
$\mathbb{D}$. For $z=(z_1,\cdots,z_n)$ and
$z'=(z'_1,\cdots,z'_n)\in\mb C^n$, denote $\lt\ll z,z'\rt\rl=z_1\ov
z_1'+\cdots+z_n\ov z_n'$ and $|z|=\lt\ll z,z\rt\rl^{1/2}$. Let
$\O_{X,Y}$ be the class of all holomorphic mappings $f$ from $X$
into $Y$. For $f\in\O_{X,Y}$, define
\begin{equation}\label{2}
|\nabla|f|(z)|=\sup_{\b\in
\mathbb{C}^n,\,|\b|=1}\lt(\lim_{t\in\mathbb{R},\,t\rightarrow0^+}\frac{|f|(z+t\b)-|f|(z)}{t}\rt),\
\ \ \ \ \ z\in X,
\end{equation}
where $f=(f_1,\cdots,f_m)$ and
$|f|=(|f_1|^2+\cdots+|f_m|^2)^{\frac{1}{2}}$. Some calculation for
$|\nabla|f||$ will be given in Section 2.

For $f\in\O_{\mathbb{D},\mathbb{D}}$, the classical Schwarz-Pick
lemma says that
\begin{equation}\label{3}
|f'(z)|\leq\frac{1-|f(z)|^2}{1-|z|^2},\ \ \ \ \ \ z\in\mathbb{D}.
\end{equation}
This inequality does not hold for $f\in\O_{\mathbb{D},\mathbb{B}_m}$
with $m\geq2$. For instance, the mapping $f(z)=(z,1)/\sqrt{2}$
satisfies
$$|f'(0)|=\sqrt{1-|f(0)|^2}>1-|f(0)|^2.$$ Note that for $f\in\O_{\mathbb{D},\mathbb{B}_m}$ with
$m\geq2$, \cite{D1} and \cite{P} proved the following inequality
\begin{equation*}
|f'(z)|\leq\frac{\sqrt{1-|f(z)|^2}}{1-|z|^2},\ \ \ \ \ \
z\in\mathbb{D}.
\end{equation*}
However Pavlovi$\acute{c}$ \cite{P} found that \eqref{3} can also be
written as
\begin{equation}\label{4}
|\nabla|f|(z)|\leq\frac{1-|f(z)|^2}{1-|z|^2},\ \ \ \ \ \
z\in\mathbb{D},
\end{equation}
since for $f\in\O_{\mathbb{D},\mathbb{D}}$, $|\nabla|f|(z)|=|f'(z)|$
by (2.12). In \cite{P}, Pavlovi$\acute{c}$ proved that this form
\eqref{4} of Schwarz-Pick lemma can be extended to
$\O_{\mathbb{D},\mathbb{B}_m}$. Precisely: for
$f\in\O_{\mathbb{D},\mathbb{B}_m}$, the following inequality holds
\begin{equation}\label{50}
|\nabla|f|(z)|\leq\frac{1-|f(z)|^2}{1-|z|^2},\ \ \ \ \ \
z\in\mathbb{D}.
\end{equation}
So it is natrual for us to consider that if the form \eqref{4} of
Schwarz-Pick lemma can also be extended to
$\O_{\mathbb{B}_n,\mathbb{B}_m}$.

In this paper, we generalize the form \eqref{4} of Schwarz-Pick
lemma to $\O_{\mathbb{B}_n,\mathbb{B}_m}$ and obtain the following
theorem: \bt\label{th1} Let $f: \mathbb{B}_n\rightarrow\mathbb{B}_m$
be a holomorphic mapping. Then
\begin{equation}\label{33}
|\nabla|f|(z)|\leq\frac{1-|f(z)|^2}{1-|z|^2},\ \ \ \ \ \
z\in\mathbb{B}_n.
\end{equation}
\et

In addition, it is well known that if equality holds in \eqref{3} at
some point $p\in\mathbb{D}$ then
\begin{equation}\label{51}
f(z)=\frac{a+e^{i\theta}z}{1+\bar{a}e^{i\theta}z}
\end{equation}
for some $a\in\mathbb{D}$, some $\theta\in\mathbb{R}$ and all
$z\in\mathbb{D}$. In this paper, we also discuss the equality case
in \eqref{33} and obtain the following theorems, which show that
\eqref{33} is sharp.

\bt\label{th3} Let $f: \mathbb{D}\rightarrow\mathbb{B}_m$ be a
holomorphic mapping. Let $p$ be a given point in $\mathbb{D}$. If
$|\nabla|f|(p)|=\frac{1-|f(p)|^2}{1-|p|^2}$, then
\begin{itemize}
\item[(1)] if $f(p)=0$,
then $$f(z)=\b\vp_p(z)$$ for some $\b\in\mathbb{C}^m$ with $|\b|=1$
and all $z\in\mathbb{D}$;
\item[(2)] if $f(p)\neq0$,
then $$\lt\ll
f(z),\frac{a}{|a|}\rt\rl=\frac{|a|+e^{i\theta}\vp_p(z)}{1+|a|e^{i\theta}\vp_p(z)}$$
for some $a\in\mathbb{B}_m$ with $a\neq0$, some
$\theta\in\mathbb{R}$ and all $z\in\mathbb{D}$;
\end{itemize}
where $\vp_p(z)=\frac{p-z}{1-\bar{p}z}$ for all $z\in\mathbb{D}$.
\et

\bt\label{th4} Let $f: \mathbb{B}_n\rightarrow\mathbb{B}_m$ be a
holomorphic mapping. Let $p$ be a given point in $\mathbb{B}_n$. If
$|\nabla|f|(p)|=\frac{1-|f(p)|^2}{1-|p|^2}$, then
\begin{itemize}
\item[(1)] if $f(p)=0$,
then
$$f(p+z(q-p))=\b\vp_{\frac{-c_{p,q}}{r_{p,q}}}\lt(\frac{z-c_{p,q}}{r_{p,q}}\rt)$$
for some $q\in\mathbb{B}_n$ with $q\neq p$ and that $q-p$ and $p$
are collinear, some $\b\in\mathbb{C}^m$ with $|\b|=1$ and all
$z\in\mathbb{D}_{c_{p,q},r_{p,q}}$;
\item[(2)] if $f(p)\neq0$,
then $$\lt\ll f(p+z(q-p)),\frac{a}{|a|}\rt\rl=
\frac{|a|+e^{i\theta}\vp_{\frac{-c_{p,q}}{r_{p,q}}}\lt(\frac{z-c_{p,q}}{r_{p,q}}\rt)}
{1+|a|e^{i\theta}\vp_{\frac{-c_{p,q}}{r_{p,q}}}\lt(\frac{z-c_{p,q}}{r_{p,q}}\rt)}$$
for some $q\in\mathbb{B}_n$ with $q\neq p$ and that $q-p$ and $p$
are collinear, some $a\in\mathbb{B}_m$ with $a\neq0$, some
$\theta\in\mathbb{R}$ and all $z\in\mathbb{D}_{c_{p,q},r_{p,q}}$;
\end{itemize}
where $c_{p,q}=-\frac{\langle p,q-p\rangle}{|q-p|^2}$,
$r_{p,q}=\sqrt{\frac{1-|p|^2}{|q-p|^2}+\left|\frac{\langle
p,q-p\rangle}{|q-p|^2}\right|^2}$,
$\mb{D}_{c_{p,q},r_{p,q}}=\{z\in\mathbb{C}: |z-c_{p,q}|<r_{p,q}\}$
 and
$\vp_{z_0}(z)=\frac{z_0-z}{1-\overline{z_0}z}$ for any
$z_0\in\mathbb{D}$ and all $z\in\mathbb{D}$. \et Note that Theorem
\ref{th4} is coincident with Theorem \ref{th3} when $n=1$.

In Section 2, some calculation for $|\nabla|f||$ will be given. In
Section 3, we will give the proof of Theorem \ref{th1}. In Section
4, we will give the proofs of Theorem \ref{th3} and Theorem
\ref{th4}.

\section{Some calculation for
$|\nabla|f||$}

For $f\in\O_{X,Y}$, let $g=|f|$. By \eqref{2}, we know that
\begin{equation}\label{5}
|\nabla g(z)|=\sup_{\b\in
\mathbb{C}^n,\,|\b|=1}\lt(\lim_{t\in\mathbb{R},\,t\rightarrow0^+}\frac{g(z+t\b)-g(z)}{t}\rt),\
\ \ \ \ \ z\in X.
\end{equation}

If $g(z)\neq0$, then $g$ is $\mathbb{R}$-differentiable at $z$ and
$\nabla g$ is the ordinary gradient. Let $z=(z_1,\cdots,z_n)$ and
$z_j=x_j+iy_j$ for $j=1,\cdots,n$, where $x_j\in\mathbb{R}$ and
$y_j\in\mathbb{R}$. Then
\begin{equation*}
\frac{\p g}{\p z_j}=\frac{1}{2}\lt(\frac{\p g}{\p x_j}-i\frac{\p
g}{\p y_j}\rt)\quad\mbox{for $j=1,\cdots,n$},
\end{equation*}
and
\begin{equation}\label{6}
\begin{split}
|\nabla g|&=\lt|\lt(\frac{\p g}{\p x_1},\frac{\p g}{\p
y_1},\cdots,\frac{\p g}{\p x_n},\frac{\p g}{\p y_n}\rt)\rt|\\
&=2\lt|\lt(\frac{\p g}{\p z_1},\cdots,\frac{\p g}{\p z_n}\rt)\rt|.
\end{split}
\end{equation}
Note that for $j=1,\cdots,n$,
\begin{equation}\label{7}
\begin{split}
\frac{\p g}{\p z_j}&=\frac{\p |f|}{\p z_j}\\
&=\frac{1}{2|f|}\frac{\p |f|^2}{\p z_j}\\
&=\frac{1}{2|f|}\frac{\p (|f_1|^2+\cdots+|f_m|^2)}{\p z_j}\\
&=\frac{1}{2|f|}\frac{\p (f_1\bar{f_1}+\cdots+f_m\bar{f_m})}{\p
z_j}\\
&=\frac{1}{2|f|}\lt(\frac{\p f_1}{\p z_j}\bar{f_1}+\cdots+\frac{\p
f_m}{\p z_j}\bar{f_m}\rt)\\
&=\frac{1}{2|f|}\lt\ll \frac{\p f}{\p z_j}, f\rt\rl,
\end{split}
\end{equation}
where $f=(f_1,\cdots,f_m)$ and $\frac{\p f}{\p z_j}=\lt(\frac{\p
f_1}{\p z_j},\cdots,\frac{\p f_m}{\p z_j}\rt)$. Then by \eqref{6}
and \eqref{7}, we have that for $g(z)\neq0$,
\begin{equation}\label{8}
|\nabla g(z)|=\frac{1}{|f(z)|}\lt|\lt(\lt\ll \frac{\p f(z)}{\p z_1},
f(z)\rt\rl,\cdots,\lt\ll \frac{\p f(z)}{\p z_n}, f(z)\rt\rl\rt)\rt|.
\end{equation}

If $g(z)=0$, then $f(z)=0$ and by \eqref{5} we have
\begin{equation}\label{9}
\begin{split}
|\nabla g(z)|&=\sup_{\b\in
\mathbb{C}^n,\,|\b|=1}\lt(\lim_{t\in\mathbb{R},\,t\rightarrow0^+}\frac{g(z+t\b)}{t}\rt)\\
&=\sup_{\b\in
\mathbb{C}^n,\,|\b|=1}\lt(\lim_{t\in\mathbb{R},\,t\rightarrow0^+}\frac{|f(z+t\b)-f(z)|}{t}\rt).
\end{split}
\end{equation}
For $f\in\O_{X,Y}$ and any $z\in X$, there is a bounded linear
operator $Df(z)$ of $\mathbb{C}^n$ into $\mathbb{C}^m$ such that
\begin{equation}\label{1}
\lim_{\b\in
\mathbb{C}^n,\,|\b|\rightarrow0}\frac{|f(z+\b)-f(z)-Df(z)\cdot \b|}
{|\b|}=0,
\end{equation}
where $Df(z)\cdot \b$ denotes the evaluation of $Df(z)$ on $\b\in
\mathbb{C}^n$. $Df(z)$ is called the Fr\'echet derivative of $f$ at
$z$, and $Df(z)\cdot \b$ is called the Fr\'echet derivative of $f$
at $z$ in the direction $\b$. In fact, for $f\in\O_{X,Y}$, $z\in X$
and $\b\in \mathbb{C}^n$, the following equality holds:
\begin{equation}\label{12}
Df(z)\cdot \b=\sum^n_{j=1}\b_j\frac{\p f(z)}{\p z_j},
\end{equation}
where $\b=(\b_1,\cdots,\b_n)$, $f=(f_1,\cdots,f_m)$ and $\frac{\p
f}{\p z_j}=\lt(\frac{\p f_1}{\p z_j},\cdots,\frac{\p f_m}{\p
z_j}\rt)$ for $j=1,\cdots,n$. Then by \eqref{1}, we have for any
$\b\in \mathbb{C}^n$ and $|\b|=1$,
\begin{equation*}
\lim_{t\in\mathbb{R},\,t\rightarrow0^+}\frac{|f(z+t\b)-f(z)-Df(z)\cdot
\b|}{t}=0.
\end{equation*}
So for any $\b\in \mathbb{C}^n$ and $|\b|=1$,
\begin{equation}\label{10}
\lim_{t\in\mathbb{R},\,t\rightarrow0^+}\frac{|f(z+t\b)-f(z)|}{t}=|Df(z)\cdot
\b|.
\end{equation}
Then by \eqref{9} and \eqref{10}, we obtain for $g(z)=0$,
\begin{equation}\label{11}
|\nabla g(z)|=\sup_{\b\in \mathbb{C}^n,\,|\b|=1}|Df(z)\cdot \b|.
\end{equation}

Now by \eqref{8} and \eqref{11}, for $f\in\O_{X,Y}$, we have
\begin{equation}\label{14}
|\nabla|f|(z)|=
\begin{cases}
\frac{1}{|f(z)|}\lt|\lt(\lt\ll \frac{\p f(z)}{\p z_1},
f(z)\rt\rl,\cdots,\lt\ll \frac{\p f(z)}{\p z_n}, f(z)\rt\rl\rt)\rt|, &\mbox{if}\quad f(z)\neq0;\\
\underset{\b\in \mathbb{C}^n,\,|\b|=1}{\sup}|Df(z)\cdot \b|,
&\mbox{if}\quad f(z)=0.
\end{cases}
\end{equation}
Then for $f\in\O_{X,Y}$ with $X\subset\mathbb{C}$,
\begin{equation}\label{15}
|\nabla|f|(z)|=
\begin{cases}
\frac{1}{|f(z)|}\lt|\lt\ll f'(z),
f(z)\rt\rl\rt|, &\mbox{if}\quad f(z)\neq0;\\
|f'(z)|, &\mbox{if}\quad f(z)=0,
\end{cases}
\end{equation}
since by \eqref{12},
$$\sup_{\b\in \mathbb{C},\,|\b|=1}|Df(z)\cdot \b|=\sup_{\b\in \mathbb{C},\,|\b|=1}|\b f'(z)|=|f'(z)|.$$
In particular, for $f\in\O_{X,Y}$ with $X\subset\mathbb{C}$ and
$Y\subset\mathbb{C}$,
\begin{equation}\label{13}
|\nabla|f|(z)|=|f'(z)|.
\end{equation}

\section{Proof of Theorem \ref{th1}}

First we give two lemmas.

\bl\label{l1}  \cite{D2} For given $p,q\in\mathbb{B}_n$ with $q\neq
p$, let $L(z)=p+z(q-p)$ for $z\in\mathbb{C}$. Then,
$$L(\mb{D}_{c_{p,q},r_{p,q}})\subset \mathbb{B}_n\ \ \ \ \
L(\partial\mb{D}_{c_{p,q},r_{p,q}})\subset
\partial\mathbb{B}_n,$$ where
\begin{equation}\label{20}
\begin{split}
\mb{D}_{c_{p,q},r_{p,q}}&=\{z\in\mathbb{C}:
|z-c_{p,q}|<r_{p,q}\},\\
c_{p,q}&=-\frac{\langle p,q-p\rangle}{|q-p|^2},\quad
r_{p,q}=\sqrt{\frac{1-|p|^2}{|q-p|^2}+\left|\frac{\langle
p,q-p\rangle}{|q-p|^2}\right|^2}.
\end{split}
\end{equation}
\el

\bl\label{l2} Let $g:\mb{D} _{c,r}\longrightarrow \mathbb{B}_m$ be a
holomorphic mapping, where
$\mb{D}_{c,r}=\{z\in\mathbb{C}:|z-c|<r\}$. Then
$$|\nabla|g|(z)|\leq\frac{r(1-|g(z)|^2)}{r^2-|z-c|^2},\ \ \ \ \ \
z\in\mb{D}_{c,r}.$$ \el

\bp Let
\begin{equation}\label{16}
h(z)=g(rz+c),\ \ \ \ \ \ z\in\mb{D}.
\end{equation}
Then $h: \mb{D}\rightarrow\mathbb{B}_m$ is a holomorphic mapping.
From \eqref{50}, we know that
$$|\nabla|h|(z)|\leq\frac{1-|h(z)|^2}{1-|z|^2},\ \ \ \ \ \
z\in\mb{D}.$$ By the above inequality and \eqref{15}, we have
\begin{equation}\label{17}
\begin{cases}
\frac{1}{|h(z)|}\lt|\lt\ll h'(z),h(z)\rt\rl\rt|\leq\frac{1-|h(z)|^2}{1-|z|^2}, &\mbox{if}\quad h(z)\neq0;\\
|h'(z)|\leq\frac{1}{1-|z|^2}, &\mbox{if}\quad h(z)=0.
\end{cases}
\end{equation}
Let $\xi=rz+c$. Then by \eqref{16} and \eqref{17}, we have
$h(z)=g(\xi)$, $h'(z)=rg'(\xi)$, and
\begin{equation}\label{18}
\begin{cases}
\frac{1}{|g(\xi)|}\lt|\lt\ll g'(\xi),g(\xi)\rt\rl\rt|\leq\frac{r(1-|g(\xi)|^2)}{r^2-|\xi-c|^2}, &\mbox{if}\quad g(\xi)\neq0;\\
|g'(\xi)|\leq\frac{r}{r^2-|\xi-c|^2}, &\mbox{if}\quad g(\xi)=0.
\end{cases}
\end{equation}
Note that by \eqref{15},
\begin{equation}\label{19}
|\nabla|g|(\xi)|=
\begin{cases}
\frac{1}{|g(\xi)|}\lt|\lt\ll g'(\xi),g(\xi)\rt\rl\rt|, &\mbox{if}\quad g(\xi)\neq0;\\
|g'(\xi)|, &\mbox{if}\quad g(\xi)=0.
\end{cases}
\end{equation}
Therefore by \eqref{18} and \eqref{19}, we obtain that
$$|\nabla|g|(\xi)|\leq\frac{r(1-|g(\xi)|^2)}{r^2-|\xi-c|^2},\ \ \ \ \ \
\xi\in\mb{D}_{c,r}.$$ The lemma is proved. \ep

Now we give the proof of Theorem \ref{th1}.

\noindent{\it Proof of Theorem \ref{th1}.}\quad  For given
$p,q\in\mathbb{B}_n$ with $q\neq p$, let
$$g(z)=f(p+z(q-p)),\ \ \ \ \ \
z\in\mb{D}_{c_{p,q},r_{p,q}},$$ where $\mb{D}_{c_{p,q},r_{p,q}}$ is
defined in Lemma \ref{l1}. Then $g:
\mb{D}_{c_{p,q},r_{p,q}}\rightarrow\mathbb{B}_m$ is a holomorphic
mapping and
\begin{equation}\label{25}
g(0)=f(p), \ \ \ \ \ g'(0)=Df(p)\cdot (q-p).
\end{equation}
By \eqref{15}, we have
\begin{equation}\label{21}
|\nabla|g|(0)|=
\begin{cases}
\frac{1}{|g(0)|}\lt|\lt\ll g'(0),g(0)\rt\rl\rt|, &\mbox{if}\quad g(0)\neq0;\\
|g'(0)|, &\mbox{if}\quad g(0)=0.
\end{cases}
\end{equation}
By Lemma \ref{l2}, we have
\begin{equation}\label{22}
|\nabla|g|(0)|\leq\frac{r_{p,q}(1-|g(0)|^2)}{r_{p,q}^2-|c_{p,q}|^2}.
\end{equation}
Note that
\begin{equation}\label{23}
\begin{split}
\frac{r_{p,q}}{r_{p,q}^2-|c_{p,q}|^2}&=\frac{1}{1-|p|^2}\lt(|q-p|
^2(1-|p|^2)+|\lt\ll p, q-p\rt\rl|^2\rt)^{\frac{1}{2}}\\
&\leq\frac{1}{1-|p|^2}|q-p|
\end{split}
\end{equation}
since \eqref{20}. Then from \eqref{22} and \eqref{23}, we have
\begin{equation}\label{24}
|\nabla|g|(0)|\leq\frac{1-|g(0)|^2}{1-|p|^2}|q-p|.
\end{equation}
By \eqref{21}, \eqref{24} and \eqref{25}, we have
\begin{equation*}
\begin{cases}
\frac{1}{|f(p)|}\lt|\lt\ll Df(p)\cdot
(q-p),f(p)\rt\rl\rt|\leq\frac{1-|f(p)|^2}{1-|p|^2}|q-p|,
&\mbox{if}\quad f(p)\neq0;\\
|Df(p)\cdot (q-p)|\leq\frac{1}{1-|p|^2}|q-p|, &\mbox{if}\quad
f(p)=0.
\end{cases}
\end{equation*}
Then
\begin{equation*}
\begin{cases}
\frac{1}{|f(p)|}\lt|\lt\ll Df(p)\cdot
\frac{q-p}{|q-p|},f(p)\rt\rl\rt|\leq\frac{1-|f(p)|^2}{1-|p|^2},
&\mbox{if}\quad f(p)\neq0;\\
|Df(p)\cdot \frac{q-p}{|q-p|}|\leq\frac{1}{1-|p|^2}, &\mbox{if}\quad
f(p)=0.
\end{cases}
\end{equation*}
So for any $\b\in\mathbb{C}^n$ with $|\b|=1$, we obtain that
\begin{eqnarray}
\begin{cases}
\frac{1}{|f(p)|}\lt|\lt\ll Df(p)\cdot
\b,f(p)\rt\rl\rt|\leq\frac{1-|f(p)|^2}{1-|p|^2},
&\mbox{if}\quad f(p)\neq0;\label{26}\\
|Df(p)\cdot \b|\leq\frac{1}{1-|p|^2}, &\mbox{if}\quad
f(p)=0.\label{27}
\end{cases}
\end{eqnarray}
For \eqref{26}, by \eqref{12} we have
\begin{equation}\label{28}
\begin{split}
\lt\ll Df(p)\cdot \b,f(p)\rt\rl&=\lt\ll \sum^n_{j=1}\b_j\frac{\p
f(p)}{\p z_j},f(p)\rt\rl\\
&=\sum^n_{j=1}\b_j\lt\ll\frac{\p f(p)}{\p
z_j},f(p)\rt\rl\\&=\sum^n_{j=1}\b_j A_j\\&=\lt\ll A, \bar{\b}\rt\rl,
\end{split}
\end{equation}
where $\b=(\b_1,\cdots,\b_n)$,
$\bar{\b}=(\bar{\b_1},\cdots,\bar{\b_n})$, $A_j=\lt\ll\frac{\p
f(p)}{\p z_j},f(p)\rt\rl$ for $j=1,\cdots,n$ and
$A=(A_1,\cdots,A_n)$. Then from \eqref{26} and \eqref{28}, we obtain
that for any $\b\in\mathbb{C}^n$ with $|\b|=1$,
\begin{equation}\label{29}
\frac{1}{|f(p)|}\lt|\lt\ll A,
\bar{\b}\rt\rl\rt|\leq\frac{1-|f(p)|^2}{1-|p|^2},\ \ \ \ \
\mbox{if}\quad f(p)\neq0.
\end{equation}
So
\begin{equation}\label{30}
\frac{1}{|f(p)|}|A|\leq\frac{1-|f(p)|^2}{1-|p|^2},\ \ \ \ \
\mbox{if}\quad f(p)\neq0,
\end{equation}
since if $|A|\neq0$, then let $\bar{\b}=\frac{A}{|A|}$ in \eqref{29}
and consequently \eqref{30} holds; if $|A|=0$, then it is obvious
that \eqref{30} holds. Then by \eqref{30} and \eqref{27}, we have
\begin{equation}\label{31}
\begin{cases}
\frac{1}{|f(p)|}|A|\leq\frac{1-|f(p)|^2}{1-|p|^2}, &\mbox{if}\quad f(p)\neq0;\\
\underset{\b\in \mathbb{C}^n,\,|\b|=1}{\sup}|Df(z)\cdot
\b|\leq\frac{1}{1-|p|^2}, &\mbox{if}\quad f(p)=0.
\end{cases}
\end{equation}
Note that by \eqref{14},
\begin{equation}\label{32}
|\nabla|f|(p)|=
\begin{cases}
\frac{1}{|f(p)|}|A|, &\mbox{if}\quad f(p)\neq0;\\
\underset{\b\in \mathbb{C}^n,\,|\b|=1}{\sup}|Df(z)\cdot \b|,
&\mbox{if}\quad f(p)=0.
\end{cases}
\end{equation}
Therefore from \eqref{31} and \eqref{32}, we have that
$$|\nabla|f|(p)|\leq\frac{1-|f(p)|^2}{1-|p|^2}.$$
The theorem is proved. \qed

\section{Proofs of Theorem \ref{th3} and Theorem \ref{th4}}

First we give the proof of Theorem \ref{th3}.

\noindent{\it Proof of Theorem \ref{th3}.}\quad First we prove the
case that $p=0$.

If $|\nabla|f|(0)|=1-|f(0)|^2$ with $f(0)=0$, then
$$|f'(0)|=1$$
since $|\nabla|f|(0)|=|f'(0)|$ by \eqref{15}. By the Lemma 1 in
\cite{D1} we have that $f(z)=f'(0)z$ for all $z\in\mathbb{D}$. So
\begin{equation}\label{36}
f(z)=\b z
\end{equation}
for some $\b\in\mathbb{C}^m$ with $|\b|=1$ and all $z\in\mathbb{D}$.
Then (1) is proved for the case that $p=0$.

If $|\nabla|f|(0)|=1-|f(0)|^2$ with $f(0)\neq0$, then by the proof
of \eqref{50} in \cite{P} we have
$$|h'(0)|=\frac{1}{|f(0)|}|\lt\ll f'(0),f(0)\rt\rl|=|\nabla|f|(0)|=1-|f(0)|^2=1-|h(0)|^2,$$
where $h(z)=\frac{1}{|f(0)|}\lt\ll f(z),f(0)\rt\rl$ is a holomorphic
function $\mathbb{D}$ from into $\mathbb{D}$. Note that
$h(0)=|f(0)|$. Then applying \eqref{51} to $h$ we get
$$h(z)=\frac{|f(0)|+e^{i\theta}z}{1+|f(0)|e^{i\theta}z}$$
for some $\theta\in\mathbb{R}$ and all $z\in\mathbb{D}$.
Consequently, $$\lt\ll
f(z),\frac{f(0)}{|f(0)|}\rt\rl=\frac{|f(0)|+e^{i\theta}z}{1+|f(0)|e^{i\theta}z}$$
for some $\theta\in\mathbb{R}$ and all $z\in\mathbb{D}$. So
\begin{equation}\label{37}
\lt\ll
f(z),\frac{a}{|a|}\rt\rl=\frac{|a|+e^{i\theta}z}{1+|a|e^{i\theta}z}
\end{equation}
for some $a\in\mathbb{B}_m$ with $a\neq0$, some
$\theta\in\mathbb{R}$ and all $z\in\mathbb{D}$. Then (2) is proved
for the case that $p=0$.

Now we prove the case that $p\neq0$. Let
\begin{equation}\label{34}
g(z)=f(\vp_p(z)),\ \ \ \ \ \ z\in\mb{D},
\end{equation}
where $\vp_p(z)=\frac{p-z}{1-\bar{p}z}$ for all $z\in\mathbb{D}$.
Then $g:\mb{D}\longrightarrow \mathbb{B}_m$ is a holomorphic mapping
and
\begin{equation}\label{35}
g(0)=f(p),\ \ \ \ g'(0)=f'(p)(-1+|p|^2).
\end{equation}

If $|\nabla|f|(p)|=\frac{1-|f(p)|^2}{1-|p|^2}$ with $f(p)=0$, then
by \eqref{15} and \eqref{35} we have $g(0)=0$ and
\begin{equation*}
\begin{split}
|\nabla|g|(0)| &=|g'(0)|\\
&=(1-|p|^2)|f'(p)|\\
&=(1-|p|^2)|\nabla|f|(p)|\\
&=1-|f(p)|^2\\
&=1-|g(0)|^2.
\end{split}
\end{equation*}
Applying \eqref{36} to $g$ we get
$$g(z)=\b z$$
for some $\b\in\mathbb{C}^m$ with $|\b|=1$ and all $z\in\mathbb{D}$.
Note that $f(z)=g(\vp_p(z))$ by \eqref{34}. Then (1) is proved for
the case that $p\neq0$.

If $|\nabla|f|(p)|=\frac{1-|f(p)|^2}{1-|p|^2}$ with $f(p)\neq0$,
then by \eqref{15} and \eqref{35} we have $g(0)\neq0$ and
\begin{equation*}
\begin{split}
|\nabla|g|(0)| &=\frac{1}{|g(0)|}|\lt\ll g'(0),g(0)\rt\rl|\\
&=\frac{1}{|f(p)|}|\lt\ll f'(p)(-1+|p|^2),f(p)\rt\rl|\\
&=(1-|p|^2)|\nabla|f|(p)|\\
&=1-|f(p)|^2\\
&=1-|g(0)|^2.
\end{split}
\end{equation*}
Applying \eqref{37} to $g$ we get
$$
\lt\ll
g(z),\frac{a}{|a|}\rt\rl=\frac{|a|+e^{i\theta}z}{1+|a|e^{i\theta}z}
$$
for some $a\in\mathbb{B}_m$ with $a\neq0$, some
$\theta\in\mathbb{R}$ and all $z\in\mathbb{D}$. Note that
$f(z)=g(\vp_p(z))$ by \eqref{34}. Then (2) is proved for the case
that $p\neq0$. The theorem is proved.
 \qed

In order to prove Theorem \ref{th4}, we need the following lemma.
\bl\label{l3} Let $g:\mb{D}_{c,r}\longrightarrow \mathbb{B}_m$ be a
holomorphic mapping, where
$\mb{D}_{c,r}=\{z\in\mathbb{C}:|z-c|<r\}$. Let $\xi$ be a given
point in $\mb{D}_{c,r}$. If
$|\nabla|g|(\xi)|=\frac{r(1-|g(\xi)|^2)}{r^2-|\xi-c|^2}$, then
\begin{itemize}
\item[(1)] if $g(\xi)=0$,
then $$g(z)=\b\vp_{\frac{\xi-c}{r}}(\frac{z-c}{r})$$ for some
$\b\in\mathbb{C}^m$ with $|\b|=1$ and all $z\in\mathbb{D}_{c,r}$;
\item[(2)] if $g(\xi)\neq0$,
then $$\lt\ll g(z),\frac{a}{|a|}\rt\rl=
\frac{|a|+e^{i\theta}\vp_{\frac{\xi-c}{r}}(\frac{z-c}{r})}{1+|a|e^{i\theta}\vp_{\frac{\xi-c}{r}}(\frac{z-c}{r})}$$
for some $a\in\mathbb{B}_m$ with $a\neq0$, some
$\theta\in\mathbb{R}$ and all $z\in\mathbb{D}_{c,r}$;
\end{itemize}
where $\vp_{z_0}(z)=\frac{z_0-z}{1-\overline{z_0}z}$ for any
$z_0\in\mathbb{D}$ and all $z\in\mathbb{D}$. \el

\bp Let
\begin{equation*}
h(z)=g(rz+c),\ \ \ \ \ \ z\in\mb{D}.
\end{equation*}
Then $h: \mb{D}\rightarrow\mathbb{B}_m$ is a holomorphic mapping and
\begin{equation}\label{38}
g(z)=h(\frac{z-c}{r}),\ \ \ \ \ \ z\in\mb{D}_{c,r}.
\end{equation}
Let $p=\frac{\xi-c}{r}$. Then $g(\xi)=h(p)$.

If $|\nabla|g|(\xi)|=\frac{r(1-|g(\xi)|^2)}{r^2-|\xi-c|^2}$ with
$g(\xi)=0$, then by the proof of Lemma \ref{l2} we have
$|\nabla|h|(p)|=\frac{1-|h(p)|^2}{1-|p|^2}$ with $h(p)=0$. Applying
Theorem \ref{th3} to $h$, we get
$$h(z)=\b\vp_p(z)$$ for some $\b\in\mathbb{C}^m$ with $|\b|=1$
and all $z\in\mathbb{D}$. Note that \eqref{38}. Then (1) is proved.

If $|\nabla|g|(\xi)|=\frac{r(1-|g(\xi)|^2)}{r^2-|\xi-c|^2}$ with
$g(\xi)\neq0$, then by the proof of Lemma \ref{l2} we have
$|\nabla|h|(p)|=\frac{1-|h(p)|^2}{1-|p|^2}$ with $h(p)\neq0$.
Applying Theorem \ref{th3} to $h$, we get
$$\lt\ll
h(z),\frac{a}{|a|}\rt\rl=\frac{|a|+e^{i\theta}\vp_p(z)}{1+|a|e^{i\theta}\vp_p(z)}$$
for some $a\in\mathbb{B}_m$ with $a\neq0$, some
$\theta\in\mathbb{R}$ and all $z\in\mathbb{D}$. Note that
\eqref{38}. Then (2) is proved. The lemma is proved. \ep

Now we give the proof of Theorem \ref{th4}.

\noindent{\it Proof of Theorem \ref{th4}.}\quad  Let
$q\in\mathbb{B}_n$ with $q\neq p$. Let
\begin{equation}\label{39}
g(z)=f(p+z(q-p)),\ \ \ \ \ \ z\in\mb{D}_{c_{p,q},r_{p,q}},
\end{equation}
where $\mb{D}_{c_{p,q},r_{p,q}}$ is defined as Lemma \ref{l1}. Then
$g: \mb{D}_{c_{p,q},r_{p,q}}\rightarrow\mathbb{B}_m$ is a
holomorphic mapping and $g(0)=f(p)$.

If $|\nabla|f|(p)|=\frac{1-|f(p)|^2}{1-|p|^2}$, then by the proof of
Theorem \ref{th1} we have $q-p$ and $p$ are collinear, and
\begin{equation}\label{40}
|\nabla|g|(0)|=\frac{r_{p,q}(1-|g(0)|^2)}{r_{p,q}^2-|c_{p,q}|^2}.
\end{equation}

If \eqref{40} holds with $g(0)=f(p)=0$, then applying Lemma \ref{l3}
to $g$, we get
$$g(z)=\b\vp_{\frac{-c_{p,q}}{r_{p,q}}}(\frac{z-c_{p,q}}{r_{p,q}})$$ for some
$\b\in\mathbb{C}^m$ with $|\b|=1$ and all
$z\in\mathbb{D}_{c_{p,q},r_{p,q}}$. Note that \eqref{39}. Then (1)
is proved.

If \eqref{40} holds with $g(0)=f(p)\neq0$, then applying Lemma
\ref{l3}  to $g$, we get
$$\lt\ll g(z),\frac{a}{|a|}\rt\rl=
\frac{|a|+e^{i\theta}\vp_{\frac{-c_{p,q}}{r_{p,q}}}(\frac{z-c_{p,q}}{r_{p,q}})}
{1+|a|e^{i\theta}\vp_{\frac{-c_{p,q}}{r_{p,q}}}(\frac{z-c_{p,q}}{r_{p,q}})}$$
for some $a\in\mathbb{B}_m$ with $a\neq0$, some
$\theta\in\mathbb{R}$ and all $z\in\mathbb{D}_{c_{p,q},r_{p,q}}$.
Note that \eqref{39}. Then (2) is proved. The theorem is proved.
 \qed

\bigbreak

{\bf Acknowledgements } This work was supported by the National
Natural Science Foundation of China (No. 11201199) and by the
Scientific Research Foundation of Jinling Institute of Technology
(No. Jit-b-201221).


\begin{thebibliography}{[99]}

\bibitem{D1} S. Y. Dai, H. H. Chen and Y. F. Pan, The Schwarz-Pick lemma of high order in several variables,
Michigan Math. J. 59(3) (2010) 517-533.

\bibitem{D2} S. Y. Dai, H. H. Chen and Y. F. Pan, The high order Schwarz-Pick lemma on complex Hilbert
balls, SCIENCE CHINA Mathematics 53(10) (2010) 2649-2656.

\bibitem{P} M. Pavlovi$\acute{c}$, A Schwarz lemma for the modulus of a vector-valued analytic
function, Proc. Amer. Math. Soc. 139(3) (2011) 969-973.

\end{thebibliography}
\end{document}